\documentclass[10pt]{article}
\usepackage{amssymb, amsmath}
\usepackage{graphicx}
\usepackage[mathscr]{eucal}

\begin{document}
\title{An analytical approach to the\\
Rational Simplex Problem}
\author{Victor Alexandrov}
\date{ }
\maketitle
\begin{abstract}
In 1973, J. Cheeger and J. Simons raised the following 
question that still remains open 
and is known as the Rational Simplex Problem:
Given a geodesic simplex in the spherical 3-space
so that all of its interior dihedral angles 
are rational multiples of $\pi$, is it true that 
its volume is a rational multiple of the volume of the 
3-sphere?
We propose an analytical approach to the 
Rational Simplex Problem by deriving a function $f(t)$,
defined as an integral of an elementary function, such 
that if there is a rational $t$, close enough to zero, 
such that the value $f(t)$ is an irrational number then 
the answer to the Rational Simplex Problem is negative.
\end{abstract}

\renewcommand{\thefootnote}{\fnsymbol{footnote}}

\section{Introduction}
\footnotetext{The author was supported in part by the 
State Maintenance Program
for Young Russian Scientists and 
the Leading Scientific Schools of the
Russian Federation (grant NSh--921.2012.1).}

At first glance, the problem 
of finding the volume of a polyhedron in the 
spherical or hyperbolic 3-space 
is both elementary and classical, and nobody can
expect new deep results here.
In fact this first impression is completely wrong.
This problem is interesting not 
only by itself  \cite{Mu12}, \cite{Vi93}, but it 
is important also 
for many problems of the geometry and topology of 3-manifolds 
\cite{AST07}, 
theory of flexible polyhedra \cite{Al97}, \cite{Sa11}, etc. 
that attract attention of modern geometers.
Several problems remain open even if we restrict our study 
by spherical 3-simplices only. One of such problems 
(that still remains open
and is known as the Rational Simplex Problem) 
was raised by Jeff Cheeger and 
James Simons in a conference held at Stanford in 1973 \cite{DS00}:
\textit{Given a geodesic simplex in the spherical 3-space
so that all of its interior dihedral angles 
are rational multiples of $\pi$, is it true that 
its volume is a rational multiple of the volume of the 
3-sphere?}

Cheeger and Simons conjectured that the answer is `usually'
no, although the answer is positive in all known cases. 
In order to demonstrate the depth of the 
Rational Simplex Problem we can mention that, 
in \cite[Theorem 2]{DS00}, Johan Dupont and Sah Chih-Han proved,
among other results, the  following theorem:
\textit{Let $\Delta$ denote a spherical simplex with all of 
its difedral angles in $\mathbb Q \pi$. Then, we have the following alternatives}:

(a) \textit{$\Delta$ is scissors congruent to a lune so  that its volume is in $\mathbb Q \pi^2$.} 

(b) \textit{$\Delta$ is not scissors congruent to a lune and exactly one of the following possibilities holds}:

{}\quad (b${}_1$) \textit{$\Delta$ leads to a negative answer of the generalized Hilbert Third Problem for spherical space;~or}

{}\quad (b${}_2$) \textit{$\Delta$ does not have volume in $\mathbb Q \pi^2$, so it is an example desired by Cheeger--Simons.}

\noindent \textit{Furthermore, for an orthoscheme $\Delta$, case} (a) \textit{and} (b) \textit{are distinguished by an explicit algorithm.} 

In this theorem, as usual, two polyhedra $P$ and $Q$ are said to be scissors
congruent if $P$ and $Q$ can each be decomposed into finite pairwase interior
disjoint unioun of polyhedra $P_i$, $Q_i$, $1\leqslant i\leqslant n$, so that
$P_i$ is congruent to $Q_i$ for each $i$.
Recall also that a \textit{lune} is defined to be the orthogonal suspension to the 
two poles of a geodesic spherical triangle lying on the equator of spherical space
and that the generalized Hilbert Third Problem still remains open and reads as follows: \textit{Is it true that two geodesic polyhedra in spherical (respectively hyperbolic) 3-space are scissors congruent if and only if they have the same volume and the same Dehn invariant?}

%
%

In the present paper we propose an analytical approach 
to the Rational Simplex Problem. 
More precisely, we introduce a function $f(t)$, 
defined as an integral of an elementary function, such 
that if there is a rational $t$, close enough to zero, 
such that the value $f(t)$ is an irrational number then 
the answer to the Rational Simplex Problem is negative.

\section{The main result}

DEFINITION. For all $|t|<1/10$, put
$$
f(t)=\frac{1}{\pi^2}\int\limits_{0}^{\pi t}
\mbox{arccos\,}\biggl(\frac{\sin s}{1+2\sin s}\biggr)ds. \eqno(1)
$$

THEOREM. \textit{If for some rational $t$, such that 
$|t|<1/10$, the number $f(t)$, given by the formula}~(1), 
\textit{is not rational then the answer to the Rational Simplex Problem 
is negative,} i.e., there exists a geodesic simplex in the 
spherical 3-space so that all of its interior dihedral angles 
are rational multiples of $\pi$, while its volume is not a 
rational multiple of $\pi^2$.

\textit{Proof:} 
Let $\mathbb S^3$ be the spherical 3-space of curvature $+1$. 
For example, we can treat $\mathbb S^3$ as the unit sphere
centered at the origin of the Euclidean space $\mathbb R^4$.

Consider a family of regular spherical simplices 
$\sigma(x)\subset\mathbb S^3$ with the edge lengths~$x$ 
varying from $0$ to $\mbox{arccos}(-1/3)$.
Observe, that

(i) for $x_*=\mbox{arccos}(-1/3)$, 
the simplex $\sigma(x_*)$ is equal to 
a hemisphere and, thus
(i${}_1$) in $\mathbb S^3$, there is no regular
spherical simplex whose edge length is greater than 
$x_*$ and
(i${}_2$) the dihedral angles $\varphi_*$ 
of the simplex $\sigma(x_*)$ are equal to~$\pi$;

(ii) the value $x_*=\mbox{arccos}(-1/3)$ is approximately 
equal to $3\pi/5$ (in fact it is approximately 1.5\% larger
than $3\pi/5$).

Equivalently, we can treat the family of regular spherical 
simplices $\sigma(x)\subset\mathbb S^3$ as being parameterized 
by their dihedral angles $\varphi$. 
In this case we write $\sigma(\varphi)$ instead of $\sigma(x)$.
Observe that $\varphi=\pi/3$ corresponds to the edge length $x=0$
and $\varphi$ increases from $\pi/3$ to $\pi$
as $x$ increases from $0$ to $x_*$.

Denote by $\alpha$ the plane angle of a face of 
$\sigma(x)$. 
A median of the face divides it into two 
mutually congruent right-angled spherical triangles; 
the spherical law of sines applied to any of them yields
$$
\frac{\sin x}{\sin(\pi/2)}=
\frac{\sin(x/2)}{\sin(\alpha/2)}
\qquad
\mbox{or}
\qquad
\cos\frac{x}{2}=
\frac{1}{2\sin(\alpha/2)}. \eqno(2)
$$
From (2), it follows that
$\alpha\to \pi/3$ as $x\to 0$ and
$\alpha\to 2\pi/3$ as $x\to x_*=\mbox{arccos}(-1/3)$.

A small sphere centered at a vertex
of  $\sigma(\varphi)$ intersects the simplex by a regular 
spherical triangle $\delta$, whose angles are equal to 
$\varphi$ radians and whose sides are equal to $\alpha$ radians.
Applying the spherical law of sines to $\delta$ in a way 
similar to the above described we get
$$
\frac{\sin\alpha}{\sin(\pi/2)}=
\frac{\sin(\alpha/2)}{\sin(\varphi/2)}
\qquad
\mbox{or}
\qquad
\cos\frac{\alpha}{2}=
\frac{1}{2\sin(\varphi/2)}. \eqno(3) 
$$
From (3), it follows that
$\varphi \to \mbox{arctan\,} 2\sqrt{2}$ as $\alpha\to \pi/3$
and
$\varphi \to \pi$ as $\alpha\to 2\pi/3$.
Observe that $\mbox{arctan\,} 2\sqrt{2}$
is, of course, the dihedral angle of a regular 
tetrahedron in the Euclidean 3-space.

Combining (2) and (3), we get 
$$
x=2\,\mbox{arccos\,}\frac{\sin(\varphi/2)}
{\sqrt{4\sin^2(\varphi/2)-1}}, \eqno(4) 
$$
where $29\pi/74\approx\mbox{arctan\,} 2\sqrt{2}<\varphi<2\pi/3$.
Observe that, in (4), $x\to 0$ as 
$\varphi\to\mbox{arctan\,} 2\sqrt{2}$
and $x\to x_*=\mbox{arccos}(-1/3)$ as $\varphi\to \pi$.
Below we will use (4) for the values of $\varphi$
that lie in the interval $(\pi/2-\pi/10, \pi/2+\pi/10)$;
this is possible according to the previous estimates.

Now recall the classical Schl{\"a}fli differential 
formula \cite{So04}:
$$
\frac{d}{dr}\mbox{vol\,}P_r= 
\frac12 \sum_j \bigl|\lambda^j_r\bigr| 
\frac{d}{dr}\beta^j_r.\eqno(5)
$$ 
Here $\{P_r\}_{r\in I}$ is a smooth
family of polyhedra in $\mathbb S^3$,
$\mbox{vol\,}P_r$ stands for the volume of $P_r$,
$\lambda^j_r$, $j=1,\dots, J$, denote 
the edges of $P_r$, 
$|\lambda^j_r|$ is the length of the edge $\lambda^j_r$, 
and  $\beta^j_r$ is the dihedral angle
of $P_r$ attached to $\lambda^j_r$.

It follows immediately from (4) and (5) that
$$
\frac{d}{d\varphi}\mbox{vol\,}\sigma(\varphi)= 
6x=
12\,\mbox{arccos\,}\frac{\sin(\varphi/2)}
{\sqrt{4\sin^2(\varphi/2)-1}}. \eqno(6)
$$

If we divide $\mathbb S^3\subset\mathbb R^4$ into
16 pieces by 4 mutually orthogonal hyperplanes in
$\mathbb R^3$, each piece will be congruent to 
$\sigma(\pi/2)$. Hence, $\mbox{vol\,}\sigma(\pi/2)=\pi^2/8$.
Taking this fact and (6) into account, we get
$$
\mbox{vol\,}\sigma(\varphi)=\frac{\pi^2}{8}+
12\int\limits_{\pi/2}^{\varphi}\mbox{arccos\,}\frac{\sin(\psi/2)}
{\sqrt{4\sin^2(\psi/2)-1}}\,d\psi \eqno(7)
$$
for all $\pi/2-\pi/10<\varphi<\pi/2+\pi/10$.
In order to simplify the integral (7), we put
$\varphi=\pi/2+\pi t$, where $|t|<1/10$, and make the 
following substitution $\psi=\pi/2+ s $. Then we get
$$
\mbox{vol\,}\sigma\biggl(\frac{\pi}{2}+\pi t\biggl)
=\frac{\pi^2}{8}+
12\int\limits_{0}^{\pi t}
\mbox{arccos\,}\frac{\sin( s /2)+
\cos( s /2)}{\sqrt{2+4\sin s }}\,d s . \eqno(8)
$$

Recall that \cite[formula 1.626.2]{RG07}:
$$
2\,\mbox{arccos\,} x=
\begin{cases}
\mbox{arccos\,} (2x^2-1), & \mbox{если\,} \hphantom{-0}0\leqslant x\leqslant 1;\\
2\pi-\mbox{arccos\,} (2x^2-1), & \mbox{если\,} -1\leqslant x\leqslant 0.
\end{cases}  \eqno(9)
$$

Taking into account that
$$
0<0.636<\frac{\sin( s /2)+
\cos( s /2)}{\sqrt{2+4\sin s }}< 0.952< 1
$$
for $-\pi/10< s <\pi/10$ and using the first line in (9),
we get from (8) after simplifications
$$
\mbox{vol\,}\sigma\biggl(\frac{\pi}{2}+\pi t\biggl)
=\frac{\pi^2}{8}+
6\int\limits_{0}^{\pi t}
\mbox{arccos\,}\frac{-\sin s }{1+2\sin s }\,d s .
$$
At last, using the formula
$\mbox{arccos\,} (-u)=\pi-\mbox{arccos\,} u$
that is valid for all $-1<u<1$, we get
$$
\mbox{vol\,}\sigma\biggl(\frac{\pi}{2}+\pi t\biggl)
=\frac{\pi^2}{8}+6\pi^2t
-6\int\limits_{0}^{\pi t}
\mbox{arccos\,}\frac{\sin s }{1+2\sin s }\,d s .
\eqno(10)
$$

It follows from (10) that if, for some rational 
$t\in(-1/10,1/10)$, the number $f(t)$, defined by 
the formula~(1),
is irrational, then the volume of the
simplex $\sigma(\pi/2+\pi t)$
is not a rational multiple of $\pi^2$.
This completes the proof of the theorem.

\section{Remarks}

\textbf{(A)} The problem to decide whether the number $f(t)$,
defined by the formula (1), is rational or not
is a difficult problem itself.
For some results on similar problems, see, e.g. \cite{Hu01,Po79, Zu05}.
The classical methods are not applicable to $f(t)$ 
at least because no representation of $f(t)$ is known
in the form of a rapidly convergent series.
By the classical methods we mean those used by 
Charles  Hermite and Ferdinand Lindemann 
on the exponential function, which allow them to prove
the transcendence of $e$ and $\pi$.
A generalization of that classical method is known as the
Siegel--Shidlovskij method for proving transcendence and 
algebraic independence of values of E-functions \cite{Sh87}.

\textbf{(B)} Note that even if, on the contrary to the 
conditions of our theorem, the number $f(t)$ is rational for all 
rational $t$ satisfying the inequality $|t|<1/10$,
we will get a new contribution to the following 
classical problem:
\textit{Given an interval $I\subset\mathbb R$, does
there exist a transcendental entire function 
$g:\mathbb C\to\mathbb C$ such that $g(x)\in\mathbb Q$
for all $x\in I\cap\mathbb Q$}?
For the first time this problem (even in a more general setting)
was solved positively by Karl Weierstrass in 1886. 
The reader may find the construction of 
Weierstrass and related historical notes in
\cite[pages 254--255]{Re91}. 
The Weierstrass's solution represents the desired 
transcendental entire function as a convergent series,
while the formula (1), of course if it is proved that 
the number $f(t)$ is rational for all 
rational $t$ satisfying the inequality $|t|<1/10$,
provides us with a closed form solution 
which is of independent interest \cite{BC13}.

\textbf{(C)}
Note also that, generally speaking, 
we may expect that what follows will lead to
a geometric proof of the rationality of the number $f(t)$,
defined by the formula (1), 
for all rational $t$ satisfying the inequality $|t|<1/10$.

By definition, let $T_0$ be the set composed of a single
simplex in the spherical 3-space
so that all its dihedral angles are rational 
multiples of $\pi$. 
Suppose that, for some $n\geqslant 0$, 
the set $T_n$ is already constructed.
Then we say that a spherical simplex $\tau$ 
belongs to the set $T_{n+1}$ if and only if
there is a spherical simplex $\tau'$ in $\cup_{k=0}^nT_k$ 
such that $\tau$ and $\tau'$
share a common face and are mutually symmetric
with respect to the 2-dimensional plane containing this face.
By definition, let $V_n$ be the set of the vertices
of the simplices $\tau\in T_n$ and let $W$ be the set of points
$x\in\mathbb S^3$ such that for every 
$\tau\in \cup_{n=0}^{\infty} T_n$ either $x\notin\tau$ or
$x$ is an interior point of $\tau$.

Observe that the following two statements 
(if proved) will provide us (among other things) with 
a geometric proof of the rationality of the number $f(t)$,
defined by the formula (1), 
for all rational $t$ satisfying the inequality $|t|<1/10$:

(h${}_1$) the set $\cup_{n=0}^{\infty} V_n$ is finite;

(h${}_2$)  for every $x,y\in W$,
the cardinality of the set of all simplices 
$\tau\in \cup_{n=0}^{\infty} T_n$
so that $x\in \tau$ is equal to 
the cardinality of the set of all simplices 
$\tau\in \cup_{n=0}^{\infty} T_n$
so that $y\in \tau$.

Of course, the hypothetical properties (h${}_1$) and (h${}_2$)
resemble the discreteness of the Coxeter groups 
\cite{Co34}, \cite{Fe04}
and the definition of the brunched covering.

\textbf{(D)} 
After reading the first version of this paper, Professor Johan Dupont
shared with me the following information that he got from Professor Ruth
Kellerhals: it turns out that the Rational Simplex Problem actually goes
back to a paper by Ludwig Schl{\"a}fli  \cite[pp. 267--269]{Sc50}, where he
has calculated the volume of several spherical orthoschemes whose
dihedral angles are rational multiples of $\pi$. Schl{\"a}fli proved that the
volume of each of those orthoschemes is a rational multiple of $\pi^2$ and
conjectured that, for all other orthoschemes, whose dihedral angles are
rational multiples of $\pi$, their volume is not a rational multiple of
$\pi^2$.

\bigskip

\noindent{Victor Alexandrov:}
\noindent\textit{Sobolev Institute of Mathematics,
Koptyug ave. 4, Novosibirsk, 630090, Russia}
\&
\textit{Department of Physics,
Novosibirsk State University, Pirogov str. 2,
Novosibirsk, 630090, Russia}

\noindent{E-mail:} \textit{alex@math.nsc.ru}

{}\hfill{Submitted: May 8, 2013}

\end{document}